\documentclass[]{elsart}
\usepackage{amsmath,amssymb%
}

\def\E{{\mathbb E}}

\def\R{{\mathbb R}}

\def\BbbN{{\mathbb N}}
\def\N{{\mathbb N}}

\def\R{\Re}

\def\bP{{\bf P}}
\def\bE{{\bf E}}

\def\bT{\mbox{\bf T}}
\def\bQ{\mbox{\bf Q}}

\def\ovf0{\widehat{f}^0}

\def\ux{\underline{x}}

\def\uX{\underline{X}}

\def\cG{{\mathcal G}}
\def\cF{{\mathcal F}}

\def\cB{{\mathcal B}}

\def\cG{{\mathcal G}}

\def\frS{{\mathfrak S}}

\def\one{{\mathbb I}}

\newtheorem{ptheorem}{Theorem}
\newtheorem{lemma}{Lemma}

\usepackage[cp1250]{inputenc}
\usepackage{fancyhdr}
\usepackage{color}
\usepackage[]{hyperref}

\begin{document}
\small
\begin{frontmatter}

\title{Optimal detection of inhomogeneous segment of observations\\
       in a stochastic sequence}
\thanks[label1]{\today}

\author{Wojciech Sarnowski\thanksref{label2}\thanksref{cor1}}
\ead{Wojciech.Sarnowski@pwr.wroc.pl}
\ead[url]{http://www.im.pwr.wroc.pl/\~{}sarnowski}
\corauth[cor1]{Corresponding author}
\author{Krzysztof Szajowski\thanksref{label2}}
\ead{Krzysztof.Szajowski@pwr.wroc.pl}
\ead[url]{http://neyman.im.pwr.wroc.pl/\~{}szajow}
\address[label2]{Wroc\l{}aw University of Technology, Institute of
Mathematics and Computer Science, Wybrze\.{z}e Wyspia\'{n}skiego 27, 50-370 Wroc\l{}aw, Poland}

\date{ \today }
\maketitle

\begin{abstract}
We register a random sequence constructed based on Markov processes by switching between them. At unobservable random
moment a change in distribution of observed sequence takes place. Using probability maximizing approach the optimal
stopping rule for detecting the disorder is identified. Some explicit solution for example is also obtained. The result is
generalization of Bojdecki's model where before and after the change independent processes are observed.
\end{abstract}


\end{frontmatter}
\fancyhf{}
\fancyhead[FC]{\thepage}
\fancyhead[HC]{}
\thispagestyle{fancy}

\section{Introduction}
This paper deals with a special problem belonging to the wide class of disorder problems. Suppose that the process $X=\{X_n,n\in\BbbN\}$, $\BbbN=\{0,1,2,\ldots\}$, is observed sequentially. It is obtained from Markov processes by switching between them at random moment $\theta$ in such a way that the process after $\theta$ starts from the state $X_{\theta-1}$. Our objective is to detect this moment based on observation of $X$. There are some papers devoted to the discrete case of such disorder detection which generalize in various directions the basic problem stated by Shiryaev~in~\cite{shi61:detection} (see e.g. Brodsky and Darkhovsky~\cite{brodar93:nonparametr}, Bojdecki~\cite{boj79:disorder}, Bojdecki and Hosza~\cite{bojhos84:problem}, Yoshida~\cite{yos83:complicated}, Szajowski~\cite{sza92:detection,sza96:twodis}).

Such model of data appears in many practical problems of the quality control~(see Brodsky and Darkhovsky~\cite{brodar93:nonparametr}, Shewhart~\cite{she31:quality} and in the collection of the papers \cite{basben86:abrupt}), traffic anomalies in networks (in papers by Dube and Mazumdar~\cite{dubmaz01:quickest}, Tartakovsky et al.~\cite{tarroz06:intrusions}), epidemiology models (see Baron~\cite{bar04:epidemio}). The aim is to recognize the moment of the change the probabilistic characteristics of the phenomenon.

Typically, disorder problem is limited to the case of switching between sequences of independent random variables (see Bojdecki \cite{boj79:disorder}). Some developments of basic model can be found in \cite{yak94:finite} where the optimal detection rule of switching moment has been obtained when the finite state-space Markov chains is disordered. Moustakides~\cite{mou98:abrupt} formulates condition which helps to reduce problem of quickest detection for dependent sequences before and after the change to the case of independent processes. Our result is generalization of results obtained by Bojdecki in \cite{boj79:disorder}. It admits Markovian dependence structure for switched sequences (with possibly uncountable state-space). We obtain an optimal rule under probability maximizing criterion.

Formulation of the problem can be found in Section \ref{sformProblem}. The main result is presented in Section \ref{rozwProblem}. Section \ref{Przyklad} provides example of application for considered model. In appendix we derive useful formulas for conditional probabilities.

\section{Formulation of the problem}\label{sformProblem}
Let $(\Omega,\mathcal{F}, \bP)$ be a probability space which supports sequence of observable random variables
$\{X_n\}_{n \in \N}$ generating filtration $\mathcal{F}_n = \sigma(X_0,X_1,...,X_n)$. The sequence takes values
in $(\E, \mathcal{B})$, where $\E$ is a subset of $\R$. Space $(\Omega,\mathcal{F}, \bP)$ supports
also unobservable (hence not measurable with respect to $\cF_n$) variable $\theta$ which has geometrical distribution:
\vspace{-3.67ex}
\begin{eqnarray}
\label{rozkladyTeta}
\bP(\theta = j) = p^{j-1}q,\;q=1-p \in (0,1),\; j=1,2,...
\end{eqnarray}
For $x \in \E$ we introduce also two homogeneous Markov processes $(X_n^0, \mathcal{G}_n^0, \bP_x^0)$,
$(X_n^1, \mathcal{G}_n^1, \bP_x^1)$ (both independent on $\theta$), which are connected with $\{X_n\}$ and $\theta$ by the following equation:
\begin{eqnarray}
\label{procesyX}
    X_n = X^0_n \cdot \one_{\{\theta>n\}} + X^1_n \cdot \one_{\{X^1_{\theta-1}=X^0_{\theta-1},\theta \leq n \}}.
\end{eqnarray}
We have that: $\cG_n^i = \sigma(X_0^i, X_1^i,\ldots ,X_n^i)$, $i \in \{0,1\}$, $n \in \{ 0, 1, 2, \ldots\}$.
On $(\E, \mathcal{B})$ for $x \in \E$ there are defined $\sigma$-additive measures $\mu(.)$ and $\mu^i_x$
($i = 0, 1$) satisfying following relations:
\begin{eqnarray*}
\bP_x^{i}(\{\omega:X_1^{i}\in B\})&=&\bP(X_1^{i}\in B|X_0^{i}=x)=\int_Bf_x^{i}(y)\mu(dy)\\
&=&\int_B\mu_x^{i}(dy)=\mu_x^{i}(B).
\end{eqnarray*}
for any $B\in\cB$.

Let us now define function $S, G$
\setlength\arraycolsep{0pt}
\begin{eqnarray}\label{gestoscLaczna}
    S(\ux_{0,n}) &=& \sum_{i=1}^n p^{i-1}qL_{n-i+1}(\ux_{0,n}) + p^{n}L_{0}(\ux_{0,n}), \\
\label{funkcjaG}
    G(\ux_{n-l-1,n},\alpha)&=& \alpha L_{l+1}(\underline{x}_{n-l-1,n}) + (1-\alpha) \\
    && \times \left( \sum_{i=0}^{l}p^{l-i}q L_{i+1}(\ux_{n-l-1,n}) + p^{l+1}L_0(\ux_{n-l-1,n})\right). \nonumber
 \end{eqnarray}
 where $x_0, x_1,\ldots, x_{n} \in \E^{n+1}, \alpha \in [0,1], 0 \leq n-l-1< n$. Here we use the following notation:
\begin{eqnarray}
    \underline{x}_{k,n} &=& (x_k, x_{k+1},...,x_{n-1},x_n),\; k\leq n, \nonumber\\
    L_{m}(\underline{x}_{k,n}) &=& \prod_{r=k+1}^{n-m} \!\!f_{x_{r-1}}^{0}(x_{r})\!\!\!\! \prod_{r=n-m+1}^{n} \!\!\!f_{x_{r-1}}^{1}(x_{r}), \nonumber\\ 
    \underline{A}_{k,n} &=&  \times_{i=k}^{n} A_{i} = A_k \times A_{k+1} \times \ldots \times A_n,\; A_i \in \cB \nonumber
\end{eqnarray}
where the convention that $\prod_{i=j_1}^{j_2}x_i = 1$ for $j_1 > j_2$ holds.

Function $S(\ux_{0,n})$ stands for join density of vector $\uX_{0,n}$. For any
$\underline{D}_{0,n}=\{\omega: \uX_{0,n} \in \underline{B}_{0,n}, B_i \in \cB\}$ and any $x\in \E$ we have:
\begin{eqnarray*}
\bP_x(\underline{D}_{0,n})=\bP(\underline{D}_{0,n}|X_0=x)&=&\int_{\underline{B}_{0,n}}S(\ux_{0,n})\mu(d\ux_{0,n})
\end{eqnarray*}
The meaning of function $G(\ux_{k,n},\alpha)$ will be clear in the sequel.

Shortly speaking our model assumes that process $\{X_n\}$ is obtained by switching at random and unknown instant $\theta$
between two Markov processes $\{X_n^0\}$ and $\{X_n^1\}$. Notice that what we assume here is that the first observation
$X_{\theta}$ after the change depends on the previous sample $X_{\theta-1}$ through the transition pdf
$f_{X_{\theta-1}}^{1}(X_{\theta})$. During on-line observation of $\{X_n\}$ we aim in detection of
switching time $\theta$ in optimal way, according to the maximum probability criterium. For any fixed
$d \in \{0,1,2,...\}$ we look for the stopping time $\tau^{*}\in \mathcal{T}$ such that
\begin{equation}
\label{PojRozregCiagowMark-Problem}
  \bP_x( | \theta - \tau^{*} | \leq d ) = \sup_{\tau \in \frS^X} \bP_x( | \theta - \tau | \leq d )
\end{equation}
where $\frS^X$ denotes the set of all stopping times with respect to the filtration
$\{\mathcal{F}_n\}_{n \in \N}$. Using parameter $d$ we control the precision level of detection.
The most rigorous case: $d=0$ will be studied in details.

\section{\label{rozwProblem}Solution of the probblem}
Let us define:
\vspace{-.5cm}
\begin{eqnarray}
Z_n &=& \bP_x(| \theta - n | \leq d \mid \mathcal{F}_n),\; n=0,1,2,\ldots, \nonumber \\
V_n &=& \rm{ess}\sup_{\{\tau \in \frS^X,\;\tau \geq n\}}\bP_x( | \theta - n | \leq d \mid \mathcal{F}_n), \; n=0,1,2,\ldots \nonumber \\
\label{stopIntuicyjny}
    \tau_0 &=& \inf\{ n: Z_n=V_n \}
\end{eqnarray}
Notice that, if $Z_{\infty}=0$, then $Z_{\tau} = \bP_x( |\theta - \tau | \leq d \mid \mathcal{F}_{\tau})$ for
$\tau \in \frS^X$. Since $\mathcal{F}_{n} \subseteq \mathcal{F}_{\tau}$ (when $n \leq \tau$) we have
\vspace{-3.67ex}
\begin{eqnarray}
    V_n &=& \rm{ess}\sup_{\tau \geq n}\bP_x(|\theta - \tau | \leq d \mid \mathcal{F}_n) = \rm{ess}\sup_{\tau \geq n}\bE_x(\one_{\{ |\theta - \tau | \leq d \}} \mid \mathcal{F}_n) \nonumber \\
            &=& \rm{ess}\sup_{\tau \geq n}\bE_x(Z_{\tau} \mid \mathcal{F}_n) \nonumber
\end{eqnarray}
The following lemma ensures existence of the solution
\begin{lemma}
    The stopping time $\tau_0$ defined by formula (\ref{stopIntuicyjny}) is the solution of problem (\ref{PojRozregCiagowMark-Problem}).
\end{lemma}

\begin{pf}
From the theorems presented in \cite{boj79:disorder} it is enough to show that
$\displaystyle{\lim_{n \rightarrow \infty}Z_n=0}$. For all natural numbers $n,k$, where $n\geq k$ we have:
\vspace{-1ex}
\begin{eqnarray*}
   Z_n &=& \bE_x(\one_{\{ |\theta - n | \leq d \}} \mid \mathcal{F}_n) \leq \bE_x(\sup_{j \geq k}\one_{\{ |\theta - j | \leq d \}} \mid \mathcal{F}_n)
\end{eqnarray*}
From Levy's theorem
$\limsup_{n\rightarrow \infty}Z_n \leq \bE_x(\sup_{j \geq k}\one_{\{ |\theta - j| \leq d \}} \mid \mathcal{F}_{\infty})$
where $\mathcal{F}_{\infty} = \sigma\left( \bigcup_{n=1}^{\infty}\mathcal{F}_n \right)$.
It is true that: $\limsup_{j \geq k,\; k\rightarrow \infty}\one_{\{ |\theta - j | \leq d \}} = 0$ \emph{a.s.}
and by the dominated convergence theorem we get
\[
  \lim_{k \rightarrow \infty}\bE_x(\sup_{j\geq k}\one_{\{ |\theta - j | \leq d \}} \mid \mathcal{F}_{\infty} ) = 0\;\; a.s.
\]
what ends the proof of the lemma.
\end{pf}

\begin{lemma}\label{PojRozregCiagowMark-lematCzasStopuMax}
Let $\tau$ be a stopping rule in the problem (\ref{PojRozregCiagowMark-Problem}). Then rule $\tilde{\tau} = \max(\tau,d+1)$
is at least as good as $\tau$. 
\end{lemma}
\begin{pf} For $\tau \geq d+1$ rules $\tau, \tilde{\tau}$ are the same. Let us consider the case when $\tau < d+1$. We
have $\tilde{\tau} = d$ and given the fact that $\bP_{x}(\theta \geq 1)=1$ we get:
\begin{eqnarray}
\bP_{x}(|\theta - \tau| \leq d) &=& \bP_{x}(\tau - d  \leq \theta \leq \tau + d)                \nonumber \\
                          &=& \bP_{x}( 1 \leq \theta \leq \tau + d)                         \nonumber \\
                          &\leq& \bP_{x}( 1 \leq \theta \leq 2d+1)                            \nonumber \\
                          &=& \bP_{x}(\tilde{\tau} - d  \leq \theta \leq \tilde{\tau} + d)  \nonumber \\
                          &=& \bP_{x}(|\theta - \tilde{\tau}| \leq d).                      \nonumber
\end{eqnarray}
\end{pf}
In consequence we can limit the class of possible stopping rules to $\frS^X_{d+1}$ i.e. stopping times equal at least
$d+1$.

For further considerations let us define posterior process:
\begin{eqnarray}
    \Pi_0 &=& 0,\nonumber\\
    \Pi_n &=& \bP_x\left(\theta \leq n \mid \mathcal{F}_n\right),\; n = 1, 2, \ldots  \nonumber
\end{eqnarray}
which is designed for information about distribution of disorder instant $\theta$. Next lemma transforms payoff function
to the more convenient form.
\begin{lemma}\label{ZmianaWyplaty}
    Let
\begin{eqnarray}\label{NowaWyplata}
h(\underline{x}_{1,d+2},\alpha) = \left( 1-p^d+ q\sum_{m=1}^{d+1}\frac{L_{m}(\underline{x}_{1,d+2}) } { p^m L_{0}(\underline{x}_{1,d+2})} \right) (1-\alpha),
\end{eqnarray}
where $x_1,...,x_{d+2} \in \E, \alpha \in (0,1)$
then
\begin{eqnarray*}
\bP_x(| \theta - n | \leq d) = \bE_x\left[ h(\underline{X}_{n-1-d,n},\Pi_n) \right]
\end{eqnarray*}
\end{lemma}

\begin{pf}
We rewrite initial criterion as the expectation
\begin{eqnarray*}
\bP_x(| \theta - n | \leq d) &=& \bE_x\left[ \bP_x(| \theta - n | \leq d \mid \mathcal{F}_n ) \right] \nonumber\\
&=& \bE_x\left[ \bP_x( \theta \leq n+d \mid \mathcal{F}_n ) - \bP_x( \theta \leq n-d-1 \mid \mathcal{F}_n ) \right]
\end{eqnarray*}
Probabilities under expectation can be transformed to the convenient form using lemmata \ref{PojRozregCiagowMark-dodatek-lemat-pin d-1 krokow wstecz}
and \ref{PojRozregCiagowMark-dodatek-lemat-pin d krokow naprzod}. Next, with the help of Lemma
\ref{PojRozregCiagowMark-dodatek-lemat-pin jako funkcja pin-d-1} (putting $l=d$) we can express
$\bP_x( \theta \leq n+d \mid \mathcal{F}_n )$ in terms of $\Pi_{n}$. Given this some straightforward calculations
imply that:
\begin{eqnarray}
\bP_{x}( |\theta - n| \leq d \mid \cF_n ) = \left( 1-p^d+q\sum_{m=1}^{d+1} \frac{L_{m}(\underline{X}_{n-d-1,n})}{p^m L_{0}(\underline{X}_{n-d-1,n})} \right)(1-\Pi_n). \nonumber
\end{eqnarray}
\end{pf}

\begin{lemma}\label{FunkcjaMarkowska}
    Process $\{\eta_n\}_{n \geq d+1}$ where $\eta_n=(\underline{X}_{n-d-1,n},\Pi_n)$ forms a random Markov function.
\end{lemma}

\begin{pf}
According to Lemma 17 pp 102-103 in \cite{shi78:optimal} it is enough to show that $\eta_{n+1}$ is a function of previous
stage $\eta_{n}$ and variable $X_{n+1}$ and that conditional distribution of $X_{n+1}$ given $\mathcal{F}_n$ is a function
of $\eta_n$. For $x_1,...,x_{d+3} \in \E, \alpha \in (0,1)$ let us consider a function

\[
    \varphi(\underline{x}_{1,d+2},\alpha,x_{d+3}) = \left(\underline{x}_{2,d+3},\frac{f_{x_{d+2}}^1(x_{d+3})(q+p\alpha)}{G(\ux_{d+2,d+3},\alpha)}\right)
\]
We will show that $\eta_{n+1} = \varphi(\eta_n,X_{n+1})$. Notice that by Lemma
\ref{PojRozregCiagowMark-dodatek-lemat-pin jako funkcja pin-d-1} ($l=0$) we get
\begin{eqnarray}\label{PojRozregCiagowMark-JedenKrokProcesPi}
\Pi_{n+1} = \frac{f_{X_n}^1(X_{n+1})(q+p\Pi_n)}{G(\uX_{n,n+1},\Pi_n)}.
\end{eqnarray}
Hence
\begin{eqnarray}
\varphi(\eta_n,X_{n+1}) &=& \varphi(\underline{X}_{n-d-1,n}, \Pi_n,X_{n+1}) \nonumber \\
&=& \left(\underline{X}_{n-d,n},X_{n+1}, \frac{f_{X_n}^1(X_{n+1})(q+p\Pi_n)}{G(\uX_{n,n+1},\Pi_n)}\right) \nonumber\\
&=& \left(\underline{X}_{n-d,n+1}, \Pi_{n+1}\right) = \eta_{n+1} \nonumber.
\end{eqnarray}
Define $\hat{\cF}_n = \sigma(\theta, \uX_{0,n})$. To see that conditional distribution of $X_{n+1}$ given $\mathcal{F}_n$ is
a function of $\eta_n$, for any Borel function $u: \E \longrightarrow \R$ let us consider the conditional expectation of
$u(X_{n+1})$ given $\mathcal{F}_n$:
\small
\begin{eqnarray}
\!\!\!\!\!\!\!\!\!\!\!\!\!\!\!\!\!\!\!\! \bE_x(&u&\!(X_{n+1})\mid \mathcal{F}_n) \nonumber\\
&=& \bE_x\left(u(X_{n+1})(1-\Pi_{n+1}) \mid \mathcal{F}_n\right) + \bE_x\left(u(X_{n+1})\Pi_{n+1} \mid \mathcal{F}_n\right) \nonumber \\
 &=& \bE_x\left(u(X_{n+1})\one_{\{\theta > n+1\}} \mid \mathcal{F}_n\right) + \bE_x\left(u(X_{n+1})\one_{\{\theta \leq n+1\}} \mid \mathcal{F}_n\right) \nonumber \\
 &=& \bE_x\left(\bE_x(u(X_{n+1})\one_{\{\theta > n+1\}}\mid \hat{\cF}_n )\mid \mathcal{F}_n\right) + \bE_x\left(\bE_x(u(X_{n+1})\one_{\{\theta \leq n+1\}}\mid \hat{\cF}_n ) \mid \mathcal{F}_n\right) \nonumber \\
&=& \bE_x\left(\one_{\{\theta > n+1\}}\bE_x(u(X_{n+1})\mid \hat{\cF}_n )\mid \mathcal{F}_n\right) + \bE_x\left(\one_{\{\theta \leq n+1\}}\bE_x(u(X_{n+1})\mid \hat{\cF}_n)  \mid \mathcal{F}_n\right) \nonumber \\
&=& \int_{\E} u(y)f^0_{X_n}(y)\mu(dy)\bP_x(\theta > n+1 \mid \cF_n) + \int_{\E} u(y)f^1_{X_n}(y)\mu(dy)\bP_x(\theta \leq n+1 \mid \cF_n) \nonumber\\
 &=& \int u(y)(p(1-\Pi_n)f^0_{X_n}(y) + (q+p\Pi_n)f^1_{X_n}(y))\mu(dy) = \int u(y)G(X_n, y, \Pi_n)\mu(dy)\nonumber
\end{eqnarray}
\normalsize
Here we use Lemma \ref{PierwszyCzlonFunkcjiWyplaty}.
\end{pf}

Lemmata \ref{ZmianaWyplaty} and \ref{FunkcjaMarkowska} are crucial for the solution of posed
problem (\ref{PojRozregCiagowMark-Problem}). They show that initial problem can be reduced to the problem of
stopping Markov random function $\eta_n = (\underline{X}_{n-d-1,n},\Pi_n )$ with the payoff given by
equation (\ref{NowaWyplata}). In consequence we can use tools of optimal stopping theory for finding stopping
time $\tau^{*}$ such that
\begin{eqnarray}\label{PrzeformulowanyProblem}
       \bE_x\left[ h(\underline{X}_{\tau^{*}-1-d,\tau^{*}},\Pi_{\tau^{*}}) \right] = \sup_{\tau \in \frS^X_{d+1}}\bE_x\left[ h(\underline{X}_{\tau-1-d,\tau},\Pi_{\tau}) \right]
\end{eqnarray}

To solve reduced problem (\ref{PrzeformulowanyProblem}) for any Borel function
$u: \E^{d+2}\times [0,1]\longrightarrow \R$ let us define operators:
\vspace{-3.67ex}
\begin{eqnarray}
  \bT u(\underline{x}_{1,d+2},\alpha) &=& \bE_{x}\left[u(\underline{X}_{n-d,n+1},\Pi_{n+1})\mid \underline{X}_{n-1-d,n}=\underline{x}_{1,d+2},\Pi_n = \alpha \right], \nonumber\\
  \bQ u(\underline{x}_{1,d+2},\alpha) &=& \max\{u(\underline{x}_{1,d+2},\alpha), \bT u(\underline{x}_{1,d+2},\alpha) \} \nonumber.
\end{eqnarray}

\begin{lemma}\label{OperatoryDlaNowejWyplaty}
    For the payoff function $h(\underline{x}_{1,d+2},\alpha)$ characterized by (\ref{NowaWyplata}) and for sequence
    $\{r_k\}_{k=0}^{\infty}$:
\small
\begin{eqnarray}
 r_0(\underline{x}_{1,d+1}) &=& p\left[ 1-p^d +q\sum_{m=1}^{d+1}\frac{L_{m-1}(\underline{x}_{1,d+1}) }{p^m L_{0}(\underline{x}_{1,d+1})}  \right], \nonumber\\
r_k(\underline{x}_{1,d+1}) &=& p\int_{\E} f_{x_{d+1}}^0(x_{d+2})\max\left\{1-p^d+q\sum_{m=1}^{d+1}\frac{L_{m}(\underline{x}_{1,d+2})}{p^m L_{0}(\underline{x}_{1,d+2})};r_{k-1}(\underline{x}_{2,d+2}) \right\}\mu(dx_{d+2}).\nonumber
\end{eqnarray}
\normalsize
the following formulas hold:
\small
\begin{eqnarray}
 \bQ^k h_1(\underline{x}_{1,d+2},\alpha) &=& (1-\alpha)\max\left\{1-p^d+q\sum_{m=1}^{d+1}\frac{L_{m}(\underline{x}_{1,d+2})}{p^m L_{0}(\underline{x}_{1,d+2})};r_{k-1}(\underline{x}_{2,d+2}) \right\},\; k \geq 1, \nonumber \\
\bT\;\bQ^k h_1(\underline{x}_{1,d+2},\alpha) &=& (1-\alpha)r_k(\underline{x}_{2,d+2}),\; k \geq 0 \nonumber.
\end{eqnarray}
\normalsize
\end{lemma}
\begin{pf}
By the definition of operator $\bT$ and using Lemma \ref{PojRozregCiagowMark-Pi_n_JakoFunkcja_Pi_n_minus_d_minus_1} ($l=0$)
given that $(\underline{X}_{n-d-1,n}, \Pi_n) = (\ux_{1,d+2},\alpha)$ we get
\footnotesize
\setlength\arraycolsep{0pt}
\begin{eqnarray}
\!\!\!\!\!\!\!\!\!\!\!\!\bT\!\!\!\!&&h(\underline{x}_{1,d+2},\alpha) = \bE_{x}\left[h(\underline{X}_{n-d,n+1},\Pi_{n+1}) \mid \underline{X}_{n-d-1,n} = \ux_{1,d+2},\Pi_n = \alpha \right] \nonumber\\
&=& \bE_{x}\left[ \bigg(1-p^d+q\sum_{m=1}^{d+1}\frac{L_{m}(\underline{X}_{n-d,n+1})}{p^m L_{0}(\underline{X}_{n-d,n+1})}\bigg)(1-\Pi_{n+1}) \mid \underline{X}_{n-d-1,n} = \ux_{1,d+2},\Pi_n = \alpha \right]\nonumber\\
&=& p(1-\alpha)\int_{\E} \left(1-p^d+q\sum_{m=1}^{d+1}\frac{L_{m-1}(\ux_{2,d+2})}{p^m L_{0}(\ux_{2,d+2})}\frac{f_{x_{d+2}}^1(x_{d+3})}{f_{x_{d+2}}^0(x_{d+3})}\right)\frac{f_{x_{d+2}}^0(x_{d+3})G(\ux_{d+2,d+3},\alpha)}{G(\ux_{d+2,d+3},\alpha)}\mu(dx_{d+3})\nonumber\\
&=&p(1-\alpha)\left[ 1-p^d+q\sum_{m=1}^{d+1}\int_{\E} \frac{L_{m-1}(\ux_{2,d+2})}{p^m L_{0}(\ux_{2,d+2})}f_{x_{d+2}}^1(x_{d+3})\mu(dx_{d+3})\right]\nonumber\\
&=&(1-\alpha)p\left[ 1-p^d+q\sum_{m=1}^{d+1}\frac{L_{m-1}(\ux_{2,d+2})}{p^m L_{0}(\ux_{2,d+2})}\right]=(1-\alpha)r_0(\underline{x}_{2,d+2}).\nonumber
\end{eqnarray}
\normalsize

Directly from the definition of $\bQ$ results that
\begin{eqnarray}
\bQ h(\underline{x}_{1,d+2},\alpha) &=& \max\left\{h(\underline{x}_{1,d+2},\alpha);\;\bT h(\underline{x}_{1,d+2},\alpha) \right\}\nonumber\\
  &=&  (1-\alpha)\max\left\{ 1-p^d+q\sum_{m=1}^{d+1}\frac{L_{m}(\underline{x}_{1,d+2})}{p^m L_{0}(\underline{x}_{1,d+2})} ;\;r_0(\underline{x}_{2,d+2})  \right\} \nonumber.
\end{eqnarray}

Suppose now that Lemma \ref{OperatoryDlaNowejWyplaty} holds for $\bT \bQ^{k-1}h$ and $\bQ^kh$ for
some $k > 1 $. Then using similar transformation as in the case of $k=0$ we get
\small
\setlength\arraycolsep{0pt}
\begin{eqnarray}
\bT &\bQ^k&h(\underline{x}_{1,d+2},\alpha) \nonumber\\
 &=& \bE_{x}\left[\bQ^kh(\underline{X}_{n-d,n+1},\Pi_{n+1}) \mid \underline{X}_{n-d-1,n} = \ux_{1,d+2},\Pi_n = \alpha\right] \nonumber\\
 &=& \int_{\E} \! \left[\max\!\! \left\{\! 1-p^d+q \!\sum_{m=1}^{d+1}\frac{L_{m}(\ux_{2,d+3})}{p^m L_{0}(\ux_{2,d+3})}; r_{k-1}(\ux_{3,d+3})\!\right\}\!(1-\alpha)pf_{x_{d+2}}^0(x_{d+3})\!\right]\!\mu(dx_{d+3})\nonumber\\ 
 &=&(1-\alpha)r_k(\ux_{2,d+2})\nonumber.
\end{eqnarray}
\normalsize
Moreover
\setlength\arraycolsep{0pt}
\begin{eqnarray}
\bQ^{k+1}h&(&\underline{x}_{1,d+2},\alpha) \nonumber\\
  &=& \max\left\{h(\underline{x}_{1,d+2},\alpha);\;\bT \bQ^kh(\underline{x}_{1,d+2},\alpha) \right\}\nonumber\\
  &=&  (1-\alpha)\max\left\{ 1-p^d+q\sum_{m=1}^{d+1}\frac{L_{m}(\underline{x}_{1,d+2})}{p^m L_{0}(\underline{x}_{1,d+2})} ;\;r_k(\underline{x}_{2,d+2})  \right\} \nonumber.
\end{eqnarray}
This completes the proof.
\end{pf}

The following theorem is the main result of the paper.
\begin{ptheorem}
\hspace{-1.5ex} \small \begin{enumerate}
  \item[(a)]
        The solution of problem (\ref{PojRozregCiagowMark-Problem}) is given by:
        \vspace{-1ex}
        \begin{eqnarray}\label{optymalnyStop}
            \tau^{*} = \inf\{n\geq d+1:1-p^d+q\sum_{m=1}^{d+1}\frac{L_{m}(\uX_{n-d-1,n})} { p^m L_{0}(\uX_{n-d-1,n})} \geq r^{*}(\underline{X}_{n-d,n})  \}
        \end{eqnarray}
        where $r^{*}(\underline{X}_{n-d,n}) = \lim_{k\longrightarrow \infty }r_k(\underline{X}_{n-d,n})$
 \item[(b)]
       Value of the problem. Given $X_0 = x$ maximal probability for (\ref{PojRozregCiagowMark-Problem}) is equal to
        \setlength\arraycolsep{0pt}
        \begin{eqnarray*}
         \bP_{x}&(&| \theta - \tau^{*} | \leq d) \nonumber\\
         &=& p^{d+1}\int_{\E^{d+1}}\max\left\{1-p^d+q\sum_{m=1}^{d+1} \frac{L_{m}(x,\underline{x}_{1,d+1})}{p^m L_{0}(x,\underline{x}_{1,d+1})}; r^{*}(\underline{x}_{1,d+1}) \right\} \nonumber\\
         && \times L_0(x,\ux_{1,d+1})\mu(d(x,\ux_{1,d+1})).
        \end{eqnarray*}
 \end{enumerate}
\end{ptheorem}

\begin{pf} Part (a).
According to Lemma \ref{PojRozregCiagowMark-lematCzasStopuMax} we look for stopping time equal at least $d+1$. From
optimal stopping theory (c.f \cite{shi78:optimal}) we know that $\tau_0$ defined by (\ref{stopIntuicyjny}) can be
expressed as
\[
    \tau_0 = \inf\{n \geq d+1: h(\underline{X}_{n-1-d,n},\Pi_n) \geq \bQ^{*}h(\underline{X}_{n-1-d,n},\Pi_n) \}
\]
where
$\bQ^{*}h(\underline{X}_{n-1-d,n},\Pi_n) = \lim_{k \longrightarrow \infty} \bQ^{k}h(\underline{X}_{n-1-d,n},\Pi_n)$.
According to Lemma \ref{OperatoryDlaNowejWyplaty}:
\setlength\arraycolsep{0.25pt}
\begin{eqnarray}
\tau_0 &=& \inf\left\{n \geq d+1: 1-p^d+q\sum_{m=1}^{d+1}\frac{L_{m}(\underline{X}_{n-d-1,n})}{p^m L_{0}(\underline{X}_{n-d-1,n})} \right. \nonumber \\
        && \;\;\;\;\;\;\;\;\;\;\;\left. \geq \max\{ 1-p^d+q\sum_{m=1}^{d+1}\frac{L_{m}(\underline{X}_{n-d-1,n})}{p^m L_{0}(\underline{X}_{n-d-1,n})} ;\;r^{*}(\underline{X}_{n-d,n})\} \right\} \nonumber\\
       &=& \inf \left\{n \geq d+1: 1-p^d+q\sum_{m=1}^{d+1}\frac{L_{m}(\underline{X}_{n-d-1,n})}{p^m L_{0}(\underline{X}_{n-d-1,n})} \geq r^{*}(\underline{X}_{n-d,n})\right\} \nonumber\\
       &=& \tau^{*}. \nonumber
\end{eqnarray}
Part (b). Basing on known facts from optimal stopping theory we can write:
{ \small
        \setlength\arraycolsep{0pt}
        \begin{eqnarray}
         \bP_{x}&(&| \theta - \tau^{*} | \leq d) \nonumber\\
         &=& \bE_{x}\left( h^{\star}_1(\underline{X}_{0,d+1},\Pi_{d+1})\right) \nonumber\\
         &=& \bE_{x}\left( (1-\Pi_{d+1})\max\left\{ 1-p^d+q\sum_{m=1}^{d+1}\frac{L_{m}(\underline{X}_{0,d+1})}{p^m L_{0}(\underline{X}_{0,d+1})} ;\;r^{\star}(\underline{X}_{1,d+1})  \right\} \right) \nonumber\\
         &=& \bE_{x}\left( \bE_{x}(\one_{\{\theta>d+1\}}\mid \cF_{d+1})\max\left\{ 1-p^d+q\sum_{m=1}^{d+1}\frac{L_{m}(\underline{X}_{0,d+1})}{p^m L_{0}(\underline{X}_{0,d+1})} ;\;r^{\star}(\underline{X}_{1,d+1})  \right\} \right) \nonumber\\
         &=& \bE_{x}\left( \one_{\{\theta>d+1\}}\max\left\{ 1-p^d+q\sum_{m=1}^{d+1}\frac{L_{m}(\underline{X}_{0,d+1})}{p^m L_{0}(\underline{X}_{0,d+1})} ;\;r^{\star}(\underline{X}_{1,d+1})  \right\}\right) \nonumber\\
         &=& \bP_{x}(\theta>d+1)\int_{\E^{d+1}}\max\left\{1-p^d+q\sum_{m=1}^{d+1} \frac{L_{m}(x,\underline{x}_{1,d+1})}{p^m L_{0}(x,\underline{x}_{1,d+1})}; r^{*}(\underline{x}_{1,d+1}) \right\} \nonumber\\
         && \times L_0(x,\ux_{1,d+1})\mu_{x_d}(d(x,\ux_{1,d+1})) \nonumber
        \end{eqnarray}
\normalsize}
What ends the proof.
\end{pf}

\section{Example}\label{Przyklad}
Let us consider the case $d=0$. Then, optimal rule (\ref{optymalnyStop}) reduces to simpler form
\begin{eqnarray*}
    \tau^{*} = \inf\{n \geq 1: \frac{f_{X_{n-1}}^1(X_n)}{pf_{X_{n-1}}^0(X_n)} \geq r^{*}(X_n) \} \nonumber
\end{eqnarray*}
with
\[
 r^{*}(X_n) = p{\int_{\E}} f_{X_{n}}^0(u) \max\{\frac{f_{X_n}^1(u)}{pf_{X_n}^0(u)}, r^{*}(u)\}d\mu(u)
\]
Moreover suppose that the state space $\E=\{0,1\}$. Matrices of transition probabilities and conditional densities
are as follow
\begin{eqnarray}
 \left[ \mu_i^0(j) \right]_{j=0,1}^{i=0,1} = \left[
            \begin{array}{cc}
                0.1\;\; & 0.9 \\
                0.8\;\; & 0.2\\
            \end{array}
    \right],\;
            &&
  \left[ \mu_i^1(j) \right]_{j=0,1}^{i=0,1} =
   \left[
            \begin{array}{cc}
                0.7\;\; & 0.3 \\
                0.4\;\; & 0.6 \\
            \end{array}
    \right] \nonumber
\end{eqnarray}

\begin{eqnarray}
 \left[ f_i^0(j) \right]_{j=0,1}^{i=0,1} = \left[
            \begin{array}{cc}
                1\;\; & 1 \\
                1\;\; & 1\\
            \end{array}
    \right],\;
            &&
  \left[ f_i^1(j) \right]_{j=0,1}^{i=0,1} =
   \left[
            \begin{array}{cc}
                7 \;\;& 1/3 \\
                1/2\;\; & 3\\
            \end{array}
    \right] \nonumber
\end{eqnarray}

For such model we find threshold $r^{*}(i)$, $i=0,1$ solving the system of equations
\[
    r^{*}(i) = \sum_{j = 0,1} pf_{i}^0(j) \max\{\frac{f_{i}^1(j)}{pf_{i}^0(j)}, r^{*}(j)\}\mu(j);\;i=0,1
\]
Treating $r^{*}$ as a function of parameter $p$ we obtain:
\begin{eqnarray}
\!\!\!\!\!\!\!\!\!\!\!\!\!r^{*}_p(0) &=& \textbf{1}_{[0,p_1]}(p) + \frac{7+9p}{10} \textbf{1}_{(p_1,p_2]}(p) + \frac{35+27p}{50-36p^2}\textbf{1}_{(p_2,p_3]}(p) + \frac{35 -7p}{50-10p-36p^2}\textbf{1}_{(p_3,1]}(p) \nonumber\\
\!\!\!\!\!\!\!\!\!\!\!\!\!r^{*}_p(1) &=& \textbf{1}_{[0,p_2]}(p) + \frac{30+28p}{50-36p^2}\textbf{1}_{(p_2,p_3]}(p) + \frac{14p}{25-50-18p^2}\textbf{1}_{(p_3,1]}(p)\nonumber
\end{eqnarray}
where:
$p_1 = \frac{1}{3}$, $p_2 = \frac{\sqrt{229}-7}{18}$, $p_3 = \frac{\sqrt{20625}-15}{136}$.
The most interesting case takes the place when $p> p_3 \approx 0,946$ because then the average disorder time is not too
small. Obtained stopping rule $\tau^{\star}$ depends on observations collected at times $\tau^{\star}-1$ and $\tau^{\star}$.
Thus, to make optimal rule more clear we need to analyze all possible sequences of $(X_{\tau^{\star}-1}, X_{\tau^{\star}})$
i.e. $\{0,0\}$, $\{0,1\}$, $\{1,0\}$, $\{1,1\}$.\\
\textbf{Sequence $\{0,0\}$:}\\
In this case we stop if only $\frac{7}{p} \geq \frac{35-7p}{50-10p-36p^2}$. Solving the inequality for $p$, we get
that stopping time takes the place for all $p \in (p_3,1)$.\\
\textbf{Sequence $\{0,1\}$:}\\
It reduces to inequality $\frac{1}{3p} \geq \frac{14p}{25-50p-18p^2}$. Taking into account that $p \in (p_3,1)$
a set of solutions is empty.\\
\textbf{Sequence $\{1,0\}$:}\\
Pair $\{1,0\}$ implies the stopping time if $\frac{7}{p} \geq \frac{35-7p}{50-10p-36p^2}$. However there is no solution
for $p \in (p_3,1)$.\\
\textbf{Sequence $\{1,1\}$:}\\
This sequence rises the alarm if only $\frac{3}{p} \geq \frac{14p}{25-50p-18p^2}$. It turns out that the inequality is
satisfied for any $p \in (p_3,1)$.

The analysis shows that we obtain very clear and simple optimal rule for case $p> p_3$:
\textbf{stop at the first moment when two "zeros" or two "ones" occur in a row}.

\appendix
\section{Lemmata}

\begin{lemma}
\label{PojRozregCiagowMark-dodatek-lemat-pin d krokow naprzod}
Let $n>0$, $k \geq 0$ then:
\begin{eqnarray}
    \label{PierwszyCzlonFunkcjiWyplaty}
\bP_{x}( \theta \leq n+k \mid \cF_n ) &=& 1 -p^k(1-\Pi_n).
\end{eqnarray}
\end{lemma}
\begin{pf} It is enough to show that for $D \in \cF_n$
\begin{eqnarray}
\int_D \one_{\{\theta >n+k\}} d\bP_x = \int_D p^k(1-\Pi_n) d\bP_x. \nonumber
\end{eqnarray}
Let us define $\widetilde{\cF}_n = \sigma(\cF_n, \one_{\{\theta >n\}})$. We have:
\setlength\arraycolsep{0.1pt}
\begin{eqnarray}
 \int_D \one_{\{\theta > n+k\}}d \bP_{x}  &=& \int_D \one_{\{\theta > n+k\}}\one_{\{\theta > n\}} d \bP_x = \int_{D \cap \{\theta > n\}} \!\!\!\one_{\{\theta > n+k\}} d \bP_x \nonumber\\
        &=& \int_{D \cap \{\theta > n\}} \!\!\!\!\!\! \bE_{x}( \one_{\{\theta > n+k\}} \mid \widetilde{\cF}_n) d \bP_x = \int_{D \cap \{\theta > n\}}\!\!\!\!\!\! \bE_{x}( \one_{\{\theta > n+k\}} \mid \theta > n) d \bP_x  \nonumber\\
        &=& \int_{D} \one_{\{\theta > n\}} p^k d \bP_x = \int_D (1-\Pi_n)p^k d \bP_x \nonumber
\end{eqnarray}
\end{pf}

\begin{lemma}
For $n > 0$ the following equality holds:
\setlength\arraycolsep{0pt}
\begin{eqnarray}
\label{PIn}
\bP_{x}&(& \theta > n \mid \cF_n ) = 1 - \Pi_n = \frac{p^nL_{0}(\underline{X}_{0,n})}{S(\underline{X}_{0,n})}.
  \end{eqnarray}
\end{lemma}
\begin{pf}
Put $\underline{D}_{0,n} = \{ \omega: \uX_{o,n} \in \underline{A}_{0,n}, A_i \in \cB\}$. Then:
\begin{eqnarray}
\bP_{x}(\underline{D}_{0,n})\bP_{x}&(&\theta>n|\underline{D}_{0,n})=\int_{\underline{D}_{0,n}}\one_{\{\theta > n\}}d\bP_x = \int_{\underline{D}_{0,n}}\bP_x(\theta > n | \cF_n)d\bP_x \nonumber \\
    &=& \int_{\underline{A}_{0,n}}\frac{p^{n} L_0(\underline{x}_{0,n})}{S(\ux_{0,n})}S(\ux_{0,n})\mu(d\ux_{0,n})= \int_{\underline{D}_{0,n}}\frac{p^{n}L_0(\underline{X}_{0,n})}{S(\uX_{0,n})}d\bP_x \nonumber
\end{eqnarray}
Hence, by definition of conditional expectation, we get the thesis.
\end{pf}

\begin{lemma}
For $\ux_{0,l+1} \in \E^{l+2}$, $\alpha \in [0,1]$ and functions $S, G$ given by equations (\ref{gestoscLaczna}) and
(\ref{funkcjaG}) we have:
 \begin{eqnarray}
 \label{FaktoryzacjaGestosci}
     S(\uX_{0,n}) &=& S(\uX_{0,n-l-1})G(\uX_{n-l-1,n},\Pi_{n-l-1})
 \end{eqnarray}
\end{lemma}
\begin{pf}
By (\ref{PIn}) we have
\small
\setlength\arraycolsep{0pt}
\begin{eqnarray}
S\!\!\!\!&(&\!\!\!\uX_{0,n-l-1})G(\uX_{n-l-1,n},\Pi_{n-l-1}) \nonumber\\
&=&S(\uX_{0,n-l-1}) \Pi_{n-l-1} L_{l+1}(\uX_{n-l-1,n})+ S(\uX_{0,n-l-1})(1-\Pi_{n-l-1})  \nonumber\\
&& \times \left( \sum_{k=0}^{l}p^{l-k}q L_{k+1}(\uX_{n-l-1,n}) + p^{l+1}L_0(\uX_{n-l-1,n})\right) \nonumber\\
&\stackrel{(\ref{PIn})}{=}& \left(\sum_{k=1}^{n-l-1}p^{k-1}q L_{n-l-k}(\uX_{0,n-l-1})\right)L_{l+1}(\uX_{n-l-1,n})+ p^{n-l-1}L_0(\underline{X}_{0,n-l-1})\nonumber\\
&& \times \left( \sum_{k=0}^{l}p^{l-k}q L_{k+1}(\uX_{n-l-1,n}) + p^{l+1}L_0(\uX_{n-l-1,n})\right)\nonumber\\
&=&\sum_{k=1}^{n-l-1}p^{k-1}q L_{n-k+1}(\uX_{0,n})+ \sum_{k=0}^{l}p^{n-k-1}q L_{k+1}(\uX_{0,n}) + p^{n}L_0(\uX_{0,n})\nonumber\\
&=&\sum_{k=1}^{n-l-1}p^{k-1}q L_{n-k+1}(\uX_{0,n})+ \sum_{k=n-l}^{n}p^{k-1}q L_{n-k+1}(\uX_{0,n}) + p^{n}L_0(\uX_{0,n})\nonumber\\
&=&\sum_{k=1}^{n}p^{k-1}q L_{n-k+1}(\uX_{0,n})+ p^{n}L_0(\uX_{0,n}) = S(\uX_{0,n}).\nonumber
\end{eqnarray}
\end{pf}

\begin{lemma}
\label{PojRozregCiagowMark-dodatek-lemat-pin d-1 krokow wstecz}
For $n >l \geq 0$ the following equation is satisfied:
\setlength\arraycolsep{0pt}
  \begin{eqnarray*}
\bP_{x}(\theta \leq n-l-1 \mid \cF_n ) &=& \frac{\Pi_{n-l-1}L_{l+1}(\underline{X}_{n-l-1,n})}{ G(\uX_{n-l-1,n},\Pi_{n-l-1})}.
  \end{eqnarray*}
\end{lemma}
\begin{pf}
Let $\underline{D}_{0,n} = \{ \omega: \uX_{o,n} \in \underline{A}_{0,n}, A_i \in \cB\}$. Then
\setlength\arraycolsep{0pt}
\begin{eqnarray*}
\bP_{x}\!\!\!&(&\!\!\!\underline{D}_{0,n})\bP_{x}(\theta>n-l-1|\underline{D}_{0,n})=\int_{\underline{D}_{0,n}}\one_{\{\theta > n-l-1\}}d\bP_x =\int_{\underline{D}_{0,n}}\bP_x(\theta > n-1 | \cF_n)d\bP_x \nonumber \\
    &=& \int_{\underline{A}_{0,n}}\frac{\sum_{k=n-l}^{n}\bP_x(\theta = k)L_{n-k+1}(\ux_{0,n}) + \bP_x(\theta > n)L_0(\ux_{0,n})}{S(\ux_{0,n})}S(\ux_{0,n})\mu(d\ux_{0,n})\nonumber \\
    &=& \int_{\underline{A}_{0,n}}\frac{p^{n-l-1} L_0(\ux_{0,n-l-1}) \left(\sum_{k=0}^{l}p^{l-k}qL_{k+1}(\ux_{n-l-1,n})+ p^{l+1}L_{0}(\ux_{n-l-1,n})\right)}{S(\ux_{0,n})} \nonumber\\
    && \times S(\ux_{0,n})\mu(d\ux_{0,n})\nonumber \\
    &=& \int_{\underline{D}_{0,n}}\frac{p^{n-l-1} L_0(\ux_{0,n-l-1}) \left(\sum_{k=0}^{l}p^{l-k}qL_{k+1}(\uX_{n-l-1,n})+ p^{l+1}L_{0}(\uX_{n-l-1,n})\right)}{S(\uX_{0,n})}d\bP_x \nonumber \\
    &\stackrel{(\ref{FaktoryzacjaGestosci})}{=}& \int_{\underline{D}_{0,n}}\frac{p^{n-l-1} L_0(\ux_{0,n-l-1}) \left(\sum_{k=0}^{l}p^{l-k}qL_{k+1}(\uX_{n-l-1,n})+  p^{l+1}L_{0}(\uX_{n-l-1,n})\right)}{S(\uX_{0,n-l-1})G(\uX_{n-l-1,n},\Pi_{n-l-1})}d\bP_x \nonumber \\
    &\stackrel{(\ref{PIn})}{=}& \int_{\underline{D}_{0,n}}(1-\Pi_{n-l-1})\frac{\sum_{k=0}^{l}p^{l-k}qL_{k+1}(\uX_{n-l-1,n})+ p^{l+1}L_{0}(\uX_{n-l-1,n})}{G(\uX_{n-l-1,n},\Pi_{n-l-1})}d\bP_x \nonumber
\end{eqnarray*}
What implies that:
\setlength\arraycolsep{0pt}
\begin{eqnarray}
\label{wzor2}
\bP_{x}&(& \theta > n-l-1 | \cF_n) \\
 &=& (1-\Pi_{n-l-1})\frac{\sum_{k=0}^{l}p^{l-k}qL_{k+1}(\uX_{n-l-1,n})+ p^{l+1}L_{0}(\uX_{n-l-1,n})}{G(\uX_{n-l-1,n},\Pi_{n-l-1})}\nonumber
\end{eqnarray}
Simple transformations of (\ref{wzor2}) lead to the thesis.
\end{pf}

\begin{lemma}
\label{PojRozregCiagowMark-dodatek-lemat-pin jako funkcja pin-d-1}
For $n >l \geq 0$ recursive equation holds:
\small
\setlength\arraycolsep{0pt}
 \begin{eqnarray}
\label{PojRozregCiagowMark-Pi_n_JakoFunkcja_Pi_n_minus_d_minus_1}
\Pi_n &=&\frac{\Pi_{n-l-1} L_{l+1}(\underline{X}_{n-l-1,n}) + (1-\Pi_{n-l-1})q\sum_{k=0}^d p^{l-k}L_{k+1}(\underline{X}_{n-l-1,n})}{ G(\uX_{n-l-1,n},\Pi_{n-l-1})}
\end{eqnarray}
\normalsize
\end{lemma}
\begin{pf}
With the aid of (\ref{PIn}) we get:
\begin{eqnarray}
\frac{1-\Pi_n}{1-\Pi_{n-l-1}} &=& \frac{p^{n}L_0(\underline{X}_{0,n})}{S(\uX_{0,n})}\frac{S(\uX_{0,{n-l-1}})}{p^{n-l-1}L_0(\underline{X}_{0,{n-l-1}})} = \frac{p^{l+1}L_{0}(\uX_{n-l-1,n})}{G(\uX_{n-l-1,n},\Pi_{n-l-1})} \nonumber
\end{eqnarray}
Hence
\begin{eqnarray}
\Pi_n &=& \frac{G(\uX_{n-l-1,n},\Pi_{n-l-1}) - p^{n-l-1}L_0(\underline{X}_{0,{n-l-1}})(1-\Pi_{n-l-1})}{G(\uX_{n-l-1,n},\Pi_{n-l-1})} \nonumber\\
    &=& \frac{\Pi_{n-l-1} L_{l+1}(\underline{X}_{n-l-1,n}) + (1-\Pi_{n-l-1})q\sum_{k=0}^d p^{l-k}L_{k+1}(\underline{X}_{n-l-1,n})}{ G(\uX_{n-l-1,n},\Pi_{n-l-1})}.\nonumber
\end{eqnarray}
\end{pf}

\normalsize
noindent\vspace{-18pt}

\end{document}